\documentclass[12pt,reqno,a4paper]{amsart}
\usepackage{amssymb}
\usepackage{color}

\usepackage[margin=1.2in]{geometry}
\def\R{\mathbb{R}}
\def\m1{{I\!\!M}}

\newcommand{\pa}{\partial}

\newcommand{\ino}{\int_{\Omega}}

\newcommand{\rife}[1]{(\ref{#1})}
\newcommand{\ov}[1]{\overline{#1}}
\newcommand{\un}[1]{\underline{#1}}

\newcommand{\txt}{\textstyle}

\renewcommand{\dfrac}{\displaystyle\frac}

\renewcommand{\i}{\infty}

\newcommand{\al}{\alpha}

\newcommand{\sg}{\sigma}

\newcommand{\om}{\Omega}
\newcommand{\lm}{\lambda}

\newtheorem{theorem}{Theorem}[section]
\newtheorem{proposition}[theorem]{Proposition}
\newtheorem{lemma}[theorem]{Lemma}
\newtheorem{corollary}[theorem]{Corollary}
\newtheorem{remark}[theorem]{Remark}
\newtheorem{definition}[theorem]{Definition}
\newcommand{\brm}{\begin{remark}\rm}
\newcommand{\erm}{\end{remark}}
\newcommand{\bdf}{\begin{definition}\rm}
\newcommand{\edf}{\end{definition}}
\newcommand{\bte}{\begin{theorem}}
\newcommand{\ete}{\end{theorem}}
\newcommand{\bpr}{\begin{proposition}}
\newcommand{\epr}{\end{proposition}}
\newcommand{\ble}{\begin{lemma}}
\newcommand{\ele}{\end{lemma}}
\newcommand{\bco}{\begin{corollary}}
\newcommand{\eco}{\end{corollary}}
\newcommand{\beq}{\begin{equation}}
\newcommand{\eeq}{\end{equation}}
\newcommand{\bdm}{\begin{displaymath}}
\newcommand{\edm}{\end{displaymath}}

\newcommand{\graf}[1]{\left\{\begin{array}{ll}#1\end{array}\right.}

\begin{document}
\numberwithin{equation}{section}
\parindent=0pt
\hfuzz=2pt
\frenchspacing

\title[Uniqueness for Toda systems]{Non-degeneracy and uniqueness of solutions to general singular Toda systems \\on bounded domains}

\author[D. Bartolucci, A. Jevnikar, J. Jin, C.S. Lin, S. Liu]{Daniele Bartolucci, Aleks Jevnikar, Jiaming Jin, Chang-Shou Lin, Senli Liu}

\address{Daniele Bartolucci, Department of Mathematics, University of Rome {\it "Tor Vergata"},  Via della Ricerca Scientifica 1, 00133 Roma, Italy.}
\email{bartoluc@mat.uniroma2.it}

\address{Aleks Jevnikar, Department of Mathematics, Computer Science and Physics, University of Udine, Via delle Scienze 206, 33100 Udine, Italy.}
\email{aleks.jevnikar@uniud.it}

\address{Jiaming Jin, School of Mathematics, Hunan University, Changsha, Hunan 410082, PR China.}
\email{jinjiaminghu@163.com}

\address{Chang-Shou Lin, Taida Institute for Mathematical Sciences and Center for Advanced Study in Theoretical Sciences, National Taiwan University, Taipei, Taiwan.}
\email{cslin@math.ntu.edu.tw}

\address{Senli Liu, School of Mathematics and Statistics, Central South University, Changsha, Hunan 410083, PR China.}
\email{mathliusl@csu.edu.cn}


\thanks{$^{(\dag)}$Research partially supported by:
Beyond Borders project 2019 (sponsored by Univ. of Rome "Tor Vergata") "{\em Variational Approaches to PDE's}",
MIUR Excellence Department Project awarded to the Department of Mathematics, Univ. of Rome Tor Vergata, CUP E83C18000100006.}

\begin{abstract}
In this note we show non-degeneracy and uniqueness results for solutions of Toda systems associated to general simple Lie algebras with multiple singular sources on bounded domains. The argument is based on spectral properties of Cartan matrices and eigenvalue analysis of linearized Liouville-type problems. This seems to be the first result for this class of problems and it covers all the Lie algebras of any rank.
\end{abstract}
\maketitle

\

\section{Introduction}

\medskip

Let $\Omega\subset\R^2$ be smooth and bounded and let $A=(a_{ij})$ be the Cartan matrix of a general simple Lie algebra. In this note we consider the following Toda system with multiple singular sources in mean field form
\begin{equation}\label{4.0}
\begin{cases}
-\Delta u_i=\sum\limits_{j=1}^na_{ij}\lambda_j\dfrac{h_je^{u_j}}{\int_{\Omega}h_je^{u_j}} \ \ &\mathrm{in} \ \ \Omega,\\
u_i=0 \ \ &\mathrm{on} \ \ \partial\Omega,
\end{cases}
\end{equation}
for $i=1,2,\cdots, n$, where $\Delta$ is the Laplace operator, $\lm_j$ are nonnegative parameters and $h_j=e^{\sg_j}$ with $\sg_j$ a singular subharmonic function having the form
\beq\label{041211.2}
\sg_i=f_i-4\pi\sum\limits_{j=1}^{m_i}\al_j^{(i)}G(x,p_j^{(i)}),
\eeq
with $f_i$ subharmonic and continuous in $\ov{\om}$, and
$\{p_1^{(i)},\cdots,p_{m_i}^{(i)}\}\subset \om$ and
$\{\al_1^{(i)},\cdots,\al_{m_i}^{(i)}\}\subset (0,+\i)$,
where
\begin{equation*}
\begin{cases}
-\Delta G(x,y)=\delta_{x=y} \ \ &\mathrm{for} \ \  x\in\Omega,\\
G(x,y)=0  \ \ &\mathrm{for} \ \  x\in\partial\Omega.
\end{cases}
\end{equation*}
Here, $\delta_p$ stands for the Dirac delta function at the point $p$.

\

All the simple Lie algebras are classified as
$$
	\textbf{A}_n, \, \textbf{B}_n, \, \textbf{C}_n, \, \textbf{D}_n,
$$
which are called classical Lie algebras and
$$
	\textbf{G}_2, \, \textbf{F}_4, \, \textbf{E}_6, \textbf{E}_7, \textbf{E}_8,
$$
which are called exceptional Lie algebras. Here, the subscript indicates the rank of the Lie algebra. For readers' convenience we explicitly write down here just the most known Cartan matrix $\textbf{A}_n$ related to $SU(n+1)$ Toda system which is given by
\begin{equation}\label{eq0}
\textbf{A}_n=\left(
\begin{array}{cccccccc}
2&-1&0&0&\cdots&0\\
-1&2&-1&0&\cdots&0\\
0&-1&2&-1&\cdots&0\\
\vdots&\vdots&\vdots&\vdots&\cdots&\vdots\\
0&\cdots&\cdots&-1&2&-1\\
0&\cdots&\cdots&0&-1&2
\end{array}\right)
\end{equation}
and refer for example to the appendix of \cite{LNW} or \cite{fh, hel, knap} for a complete list of Cartan matrices and further basic Lie theory.

\

Observe that for the Lie algebra $\textbf{A}_1$ the Cartan matrix is just $(a_{ij})=(2)$ and the Toda system reduces to the standard Liouville equation
\begin{equation}\label{liouv}
\begin{cases}
-\Delta u=2\lambda\dfrac{he^{u}}{\int_{\Omega}he^{u}} \ \ &\mathrm{in} \ \ \Omega,\\
u=0 \ \ &\mathrm{on} \ \ \partial\Omega,
\end{cases}
\end{equation}
which is related to turbulent Euler flows, abelian Chern-Simons theory \cite{J60,J61,J68} or the prescribed Gaussian curvature problem (with conic singularities), see \cite{clmp1, clmp2, kw,troy}, and has been intensively studied.

\

The Toda system also plays an important role in geometry and mathematical physics. For example, it appears in the study of holomorphic curves in projective spaces, Pl\"{u}cker formulas, harmonic maps, $\mathcal{W}$-algebras and it is a model for non-abelian Chern-Simons theory, see \cite{bfo,J9,J13,dol,fortw,gh,guest,J60,J68}.

\smallskip

For what concerns blow-up analysis, classification issues and existence results for Toda system we refer to \cite{bjmr,J34,M32,G32,Y17,Y21} and the references therein.

\

We are interested here in non-degeneracy and uniqueness of solutions to general singular Toda systems. Here, by non-degeneracy we mean that the first eigenvalue of the corresponding linearized problem is strictly positive. Up to know, this topic has been investigated only for the standard Liouville case \eqref{liouv} in various settings \cite{BJL,BL,BL2,D11, suz}. In our context, by using Alexandrov-Bol inequality and by a quite delicate eigenvalues analysis of linearized singular Liouville-type problems it was proven that if $\lm\leq4\pi$ \eqref{liouv} admits at most one solution which is non-degenerate, see in particular \cite{BL}. We also point out that this is sharp and that uniqueness does not hold in general if $\lm>4\pi$ \cite{bdm,D11}. We will not focus here on the existence of such solution since this is always granted at least in the so-called coercive regime $\lm<4\pi$ where one can just minimize the associated energy functional by means of Moser-Trudinger-type inequalities (the same being true also for the system case).

\medskip

On the other hand, it is hard to extend this kind of analysis to the Toda system and there are basically no results in this direction. The only result we are aware of is \cite{G} where a very special case with $n=2$ is treated by using the sphere covering inequality. In this note we present a simple strategy to show non-degeneracy and uniqueness of solutions for Toda systems associated to general simple Lie algebras with multiple singular sources. We stress we can cover all the Lie algebras of any rank.

\smallskip

To state our result we first need to distinguish between symmetric and non-symmetric Lie algebras. The latter are more delicate to handle but we will take advantage of a symmetric reformulation. More precisely, we will decompose the Cartan matrix $A$ of a Lie algebra as
\begin{equation} \label{dec1}
	A=DA^s,
\end{equation}  
where $D$ is a diagonal matrix and $A^s$ a symmetric matrix. Letting $d_i$ be the diagonal entries of the matrix $D$ we will also denote
\begin{equation} \label{dec2}
	\lm_i^s=d_i\lm_i, \quad i=1,\dots,n.
\end{equation}
For the symmetric Lie algebras $\textbf{A}_n, \textbf{D}_n, \textbf{E}_6, \textbf{E}_7, \textbf{E}_8$ we clearly have $A^s=A$ and $\lm_i^s=\lm_i$. The decomposition of the non-symmetric $\textbf{B}_n, \textbf{C}_n, \textbf{G}_2, \textbf{F}_4$ is postponed to subsection \ref{subsec-non-symm}.

\smallskip

We next recall that the spectral radius of a square matrix $A$ is
\begin{equation} \label{spectral}
	\rho(A)=\max\Bigr\{ |\xi| \,:\, \mbox{$\xi$ is an eigenvalue of $A$}  \Bigr\}.
\end{equation}
Observe that the known non-degeneracy and uniqueness threshold for the standard Liouville case \eqref{liouv} (Lie algebra $\textbf{A}_1$) can be written as
$$
	\lm\leq \frac{8\pi}{\rho(\textbf{A}_1)}.
$$	
We are going to extend the latter formula to any Lie algebra and obtain the following.
\begin{theorem}\label{t3}
Let $A$ be the Cartan matrix of a general simple Lie algebra and let $A^s, \lm_i^s$ be as in \eqref{dec1}, \eqref{dec2}. Suppose
\begin{equation} \label{value}
	\lm_i^s\leq \frac{8\pi}{\rho(A^s)}, \quad i=1,\dots,n.
\end{equation}	
Then, the Toda system \rife{4.0} admits at most one solution which is non-degenerate. In particular, the same holds true if
$$
	\lambda_i^s\leq2\pi, \quad i=1,\cdots,n.
$$	
\end{theorem}

\medskip

The above result holds for more general positive-definite matrices as far as some a priori estimates for solutions to the system \eqref{4.0} are available.

\medskip

Just to make an example, take now the Lie algebra  $\textbf{A}_2$ which has been studied by many authors. We have $\rho(\textbf{A}_2)=3$ and thus the threshold is given by
$$
	\lm_i\leq \frac{8\pi}{3}, \quad i=1,2.
$$
At this point let us comment about the sharpness of our result. Observe that, contrary to the standard Liouville case \eqref{liouv}, the threshold in \eqref{value} is strictly smaller than the coercivity threshold related to the Moser-Trudinger constant and one may wonder if it is possible to extend the uniqueness property in all the coercivity regime. Actually, this is not possible for example for affine Toda systems, in particular sinh-Gordon equations, see \cite{sst}, where the authors exhibit multiple solutions below the coercivity threshold. Moreover, the already mentioned results in \cite{G} (even if valid just for a special case) suggest the uniqueness threshold in \eqref{value} might be sharp, at least in some cases. This remains an interesting open problem.

\

The idea behind the proof of Theorem \ref{t3} is first to notice that a solution of the Toda system \eqref{4.0} is a subsolution of the singular Liouville-type equation
\begin{equation*}
\begin{cases}
-\Delta u_i -2K_i e^{u_i}=0&\quad \mbox{in}\;\om,\\
u_i=0&\quad\mbox{on}\;\partial\om,
\end{cases}
\end{equation*}
for $i=1,\dots,n$, where
$$
K_i={\lm}_i\dfrac{h_i}{\ino h_i e^{u_i}},
$$
for which we have a good understanding, in particular concerning eigenvalues
of its linearized problem, see for example \cite{BL} or the discussion in the sequel. We will then exploit the structure of the Toda system to derive some estimates of such eigenvalues, which are expressed in terms of the spectral radius of the associated Cartan matrix.

\smallskip

This argument will lead to the proof of non-degeneracy first. Then, the uniqueness of solutions will follow by standard arguments using the implicit function theorem and some uniform estimates for solutions to \eqref{4.0}.

\

The organization of this paper is as follows. In section \ref{sec:linear} we introduce the linearized problem and collect some useful information and in section \ref{sec:proof} we provide the proof of the non-degeneracy and uniqueness of solutions.

\

\section{{\bf The linearized problem}} \label{sec:linear}

\medskip

In this section we introduce the linearized problem which will be studied in the next section when proving non-degeneracy of solutions.

\smallskip

We start by collection some useful results about the eigenvalue analysis of linearized problems of subsolutions to singular Liouville equations. Although not stated in this generality, the following Lemmas have been proved in \cite{BL}.
\ble[\cite{BL} Theorem 2.3 and Proposition 3.2]\label{BL-1}
Let $K=e^{\sg}$ with $\sg$ taking the form \rife{041211.2} and $v$ be a smooth subsolution for
\beq\label{Liou}
\graf{
-\Delta v -K e^{\txt v}=0&\quad \mbox{in}\ \ \om,\\
v=0&\quad\mbox{on}\ \ \partial\om.
}
\eeq
If $\ino K e^{\txt v}\leq 4\pi$,
then the first eigenvalue $\nu_1$ for the linear problem
$$
\graf{
-\Delta \phi -K e^{\txt v}\phi =\nu K e^{\txt v} \phi &\quad \mbox{in}\ \ \om,\\
\phi=0&\quad\mbox{on}\ \ \partial\om,}
$$
is strictly positive.

\smallskip

Moreover, if $\ino K e^{\txt v}\leq 8\pi$,
then the second eigenvalue $\nu_2$ is strictly positive.
\ele

\smallskip

\ble[\cite{BL} Proposition 3.3]\label{lem.last}
Let $K=e^{\sg}$ with $\sg$ taking the form \rife{041211.2} and $v$ be a smooth subsolution for
\rife{Liou}. If $\ino K e^{\txt v}\leq 8\pi$ and $\phi$ solves
$$
\graf{
-\Delta \phi -K e^{\txt v}\phi =\nu K e^{\txt v} \phi &\quad \mbox{in}\ \ \om,\\
\phi=c&\quad\mbox{on}\ \ \partial\om,\\
\ino K e^{\txt v}\phi=0,}
$$
for some $c\in\R$, then $\nu$ is strictly positive.
\ele

\

We next consider the linearized problem of the Toda system
$$
\begin{cases}
-\Delta u_i=\sum\limits_{j=1}^na_{ij}\lambda_j\dfrac{h_je^{u_j}}{\int_\Omega h_je^{u_j}} \ \ &\mathrm{in} \ \ \Omega,\\
u_i=0 \ \ &\mathrm{on} \ \ \partial\Omega,
\end{cases}
$$
with $i=1,\cdots,n$, which reads
$$
\begin{cases}
-\Delta\tilde{\phi}_i=\sum\limits_{j=1}^na_{ij}\lambda_j\dfrac{h_je^{u_j}}{\int_{\Omega}h_je^{u_j}}\left(\tilde{\phi}_j-\dfrac{\int_{\Omega}h_je^{u_j}\tilde{\phi}_j}{\int_{\Omega}h_je^{u_j}}\right) \ \ &\mathrm{in} \ \ \Omega,\\
\tilde{\phi}_i=0 \ \ &\mathrm{on} \ \ \partial\Omega,
\end{cases}
$$
for $i=1,\cdots,n$.
Putting
$$
V_j=\lambda_j\frac{h_je^{u_j}}{\int_{\Omega}h_je^{u_j}},
$$
then we have
$$
\int_{\Omega}V_j=\lambda_j,
$$
and the functions
$$
\phi_j=\tilde{\phi}_j-\frac{\int_{\Omega}h_je^{u_j}\tilde{\phi}_j}{\int_{\Omega}h_je^{u_j}},
$$
satisfy
\begin{equation}\label{sys1.0}
\begin{cases}
-\Delta\phi_i=\sum\limits_{j=1}^na_{ij}V_j\phi_j, \ \ &\mathrm{in} \ \ \Omega,\\
\phi_i=c_i\in\mathbb{R}, \ \ &\mathrm{on} \ \ \partial\Omega,\\
\int_{\Omega}V_i\phi_i=0,
\end{cases}
\end{equation}
for $i=1,\cdots,n$.

\medskip

Let now $H_n=H_0^1(\Omega)\times H_0^1(\Omega)\cdots\times H_0^1(\Omega)$ and let $A^{-1}=(a^{ij})$ be the inverse matrix of $A$.
We say that $\un{\phi}=(\phi_1,\phi_2,\cdots,\phi_n)$ is a weak solution for
\rife{sys1.0} if $\un{\phi}-(c_1,c_2,\cdots,c_n)=(\phi_1,\phi_2,\cdots,\phi_n)-(c_1,c_2,\cdots,c_n)\in H_n$ satisfies
\begin{equation}\label{sys1.weak}
\ino\left[\sum\limits_{i,j=1}^n a^{ij}\nabla \phi_j\cdot\nabla \psi_i
-\sum\limits_{i=1}^n V_i\phi_i\psi_i\right]=0, \quad \forall\,\un{\psi}\in H_n,
\end{equation}
where $\un{\psi}=(\psi_1,\cdots,\psi_n)$ and
$\phi_1,\phi_2,\cdots,\phi_n$ satisfy the integral constraints in \rife{sys1.0}.

\

\section{ \textbf{The proof of Theorem \ref{t3}} } \label{sec:proof}

\medskip

In this section we derive the main Theorem \ref{t3}, starting from the symmetric Lie algebras. We divide the proof into three parts, showing first non-degeneracy and then uniqueness of solutions. We discuss about the spectral radius of Cartan matrices in the third part, proving the last assertion of Theorem \ref{t3}. The discussion of the non-symmetric case is postponed to the last subsection.

\subsection{Non-degeneracy} The proof of the non-degeneracy of solutions is performed in two steps. We consider here the symmetric case $A^s=A$ and $\lm_i^s=\lm_i$ in \eqref{dec1}, \eqref{dec2}.

\medskip

\textbf{Step 1.} We start by considering \eqref{sys1.0} with $c_i=0$ for all $i=1,\dots,n$, that is
\begin{equation}\label{linear0}
\begin{cases}
-\Delta\phi_i=\sum\limits_{j=1}^na_{ij}V_j\phi_j, \ \ &\mathrm{in} \ \ \Omega,\\
\phi_i=0, \ \ &\mathrm{on} \ \ \partial\Omega,\\
\int_{\Omega}V_i\phi_i=0,
\end{cases}
\end{equation}
for $i=1,\cdots,n$. We claim that if 
$$
\lm_i\leq \frac{4\pi}{\rho(A)}, \quad i=1,\dots,n,
$$
then, \eqref{linear0} admits only the trivial solution.

\smallskip

We prove the claim by contradiction. To avoid repetitions, we work out an argument which will be then also exploited in Step 2 when treating the general problem \eqref{sys1.0}. Consider a non-trivial weak solution of \eqref{sys1.0}, that is $H_n\ni(\phi_1,\phi_2,\cdots,\phi_n)-(c_1,c_2,\cdots,c_n)\neq(0,\cdots,0)$ satisfying \eqref{sys1.weak}. Taking
$$
(\psi_1,\cdots,\psi_n)=(\phi_1,\phi_2,\cdots,\phi_n)-(c_1,c_2,\cdots,c_n)
$$
as a test function in \eqref{sys1.weak} and using also that
$$
\int_{\Omega}V_i\phi_i=0, \quad i=1\cdots,n,
$$
we get
\begin{equation} \label{test}
\ino \left[\sum\limits_{i,j=1}^na^{ij}\nabla \phi_j\cdot\nabla \phi_i-\sum\limits_{i=1}^n V_i\phi_i^2\right] = 0.
\end{equation}
Observe now that $A$ is a positive-definite symmetric matrix and thus, recalling the spectral radius $\rho(A)$ defined in \eqref{spectral} we clearly have
$$
	\sum\limits_{i,j=1}^n\ino a^{ij}\nabla \phi_j\cdot\nabla \phi_i \geq \dfrac{1}{\rho(A)} \sum\limits_{i=1}^n\ino |\nabla\phi_i|^2
$$
and then, by \eqref{test} we get
\begin{equation} \label{negative}
	\sum\limits_{i=1}^n\left[\ino |\nabla\phi_i|^2-\rho(A) \ino V_i\phi_i^2\right] \leq 0.
\end{equation}
This argument clearly works whenever $(\phi_1,\phi_2,\cdots,\phi_n)\in H_n$ is a weak solution of \eqref{linear0}. Therefore, for any such solution, \eqref{negative} holds true. Define now
\begin{equation} \label{func}
I^{(i)}(\phi):=\int_{\Omega}|\nabla\phi|^2-\rho(A)\int_{\Omega}V_i\phi^2.
\end{equation}
The estimate in \eqref{negative} implies necessarily that
\begin{equation}\label{negative2}
I^{(i)}(\phi_i)\leq 0 \quad \mbox{for some $i\in\{1,\cdots,n\}$}.
\end{equation}
This gives an estimate about the sign of the first eigenvalue related to \eqref{func}, that is (recall $c_i=0$ in \eqref{linear0})
\beq\label{mudef}
\hat{\mu}^{(i)}_1=\inf\Bigr\{J^{(i)}\left(\phi\,\right)\,:\,\phi\in H_0^{1}(\om)\setminus\{0\}\Bigr\},
\eeq
where
$$
J^{(i)}\left(\phi\,\right)=\frac{I^{(i)}\left(\phi\,\right)}{\ino V_i \phi^2}.
$$
By \eqref{negative2} we derive that
$$
\min_{i=1,\cdots,n}\hat{\mu}^{(i)}_1\leq 0.
$$
Hence, without out loss of generality, let us assume that $\hat{\mu}^{(1)}_1\leq 0$ and let $\Phi$ be the corresponding minimizer,
which therefore satisfies
\begin{equation} \label{negative3}
\begin{split}
\begin{cases}
-\Delta \Phi -\rho(A)V_1\Phi =\hat{\mu}^{(1)}_1V_1 \Phi &\quad \mbox{in}\;\om,\\
\Phi=0&\quad\mbox{on}\;\partial\om,
\end{cases}
\end{split}
\end{equation}
This, jointly with $\lm_i\leq \dfrac{4\pi}{\rho(A)}$, will lead to a contradiction.

\medskip

Indeed, we first claim that $\rho(A)\geq 2$ for a Cartan matrix $A$ associated to any simple Lie algebra. 
This will be proved in subsection \ref{subsec-spec}.

We then observe that any solution of the Toda system \eqref{4.0} satisfies
$$
-\Delta u_i -\rho(A)K_i e^{u_i}\leq
-\Delta u_i -2K_i e^{u_i}\leq 0,\quad i=1,\cdots,n,
$$
and hence $u_i$ are smooth subsolutions of
\begin{equation}\label{subsol}
\begin{split}
\begin{cases}
-\Delta u_i -\rho(A)K_i e^{u_i}=0 \ \ &\mbox{in}\;\;\om,\\
u_i=0 \ \ &\mbox{on}\;\pa\om.
\end{cases}
\end{split}
\end{equation}
for $i=1,\cdots,n$, where
$$
K_i={\lm}_i\dfrac{h_i}{\ino h_i e^{u_i}}.
$$
Moreover, by our assumption $\lm_i\leq \dfrac{4\pi}{\rho(A)}$ we have
$$
\rho(A)\ino K_i e^{u_i}=\rho(A)\lm_i\in[0,4\pi],\quad i=1,\cdots,n.
$$
Therefore, it follows from Lemma \ref{BL-1} that the first eigenvalues $\hat{\nu}^{(i)}_1$ of
the linearized problems
\begin{equation*}
\begin{cases}
-\Delta \phi -\rho(A)V_i\phi =\hat{\nu}^{(i)}V_i \phi &\quad \mbox{in}\;\om,\\
\phi=0&\quad\mbox{on}\;\partial\om,
\end{cases}
\end{equation*}
are strictly positive, for $i=1,\cdots,n$, which is in contradiction to \eqref{negative3} with $\hat{\mu}^{(1)}_1\leq 0$.

\medskip

\textbf{Step 2.} We consider now the general linearized system
$$
\begin{cases}
-\Delta\phi_i=\sum\limits_{j=1}^na_{ij}V_j\phi_j, \ \ &\mathrm{in} \ \ \Omega,\\
\phi_i=c_i\in\mathbb{R}, \ \ &\mathrm{on} \ \ \partial\Omega,\\
\int_{\Omega}V_i\phi_i=0,
\end{cases}
$$
for $i=1,\cdots,n$, with $c_i$ non necessarily zero and claim that it admits only the trivial solution if
$$
\lm_i\leq \frac{8\pi}{\rho(A)}, \quad i=1,\dots,n.
$$

\smallskip

Indeed, reasoning as in Step 1, we have \eqref{negative2} which yields
$$
\min_{i=1,\cdots,n}\hat{\mu}^{(i)}_2\leq 0.
$$
where, with the same definition of $I^{(i)}$ and $J^{(i)}$ as above, we let
\beq\label{mudef2}
\hat{\mu}^{(i)}_2=\inf\left\{J^{(i)}\left(\phi\,\right)\,;\,\phi-c_i\in H_0^{1}(\om),\,\ino V_i\phi=0\right\}.
\eeq
Again, let us assume that $\hat{\mu}^{(1)}_2\leq 0$ and let $\Phi$ be the corresponding minimizer satisfying
$$
\graf{ -\Delta \Phi-\rho(A)V_1\Phi=\hat{\mu}^{(1)}_2V_1\Phi&\quad \mbox{in}\;\;\;\om\\
\Phi=c_1&\quad \mbox{on}\;\;\;\partial\om\\
\ino V_1\Phi=0.
}
$$
Finally, as in Step 1 this is in contradiction to Lemma \ref{lem.last} which states that $\hat{\mu}^{(i)}_2>0$ for any $i=1,\cdots,n$ since $u_i$ are smooth subsolutions of \eqref{subsol} with
$$
\rho(A)\ino K_i e^{u_i}=\rho(A)\lm_i\in[0,8\pi],\quad i=1,\cdots,n,
$$
by our assumption. Thus, we are done.

\subsection{Uniqueness} We finally prove here the solutions to the Toda system \eqref{4.0} are unique for
$$
\lm_i\leq \frac{8\pi}{\rho(A)}, \quad i=1,\dots,n,
$$
still in the symmetric case $A^s=A$ and $\lm_i^s=\lm_i$ in \eqref{dec1}, \eqref{dec2}.

\smallskip

Once we have the non-degeneracy, uniqueness of solutions follows by standard arguments so we will be sketchy, referring for example to \cite{BL} for further details. Indeed, at this point we know that the linearized operator for \eqref{4.0} has strictly positive first eigenvalue for $\lm_i\leq \frac{8\pi}{\rho(A)}$. By standard bifurcation theory, for any $\lm_i$ small enough, there exists a unique solution, bifurcating from the trivial solution $(u_i,\lm_i)=(0,0)$, $i=1,\cdots,n$. Since the linearized operator has strictly positive first eigenvalue we can apply the implicit function theorem to extend uniquely this branch of solutions for any $\lm_i\leq \frac{8\pi}{\rho(A)}$.

\smallskip

Suppose for a moment the solutions of \eqref{4.0} are uniformly bounded for $\lm_i\leq \frac{8\pi}{\rho(A)}$. Now, if by contradiction there would exist a second non-bending branch of solutions, then the uniform estimates of such solutions would force the latter branch to end up into the trivial solution $(u_i,\lm_i)=(0,0)$, $i=1,\cdots,n$, which is not possible.

\smallskip

We are thus left with the a priori estimates for solutions to \eqref{4.0} for $\lm_i\leq \frac{8\pi}{\rho(A)}$. Since the rank $1$ case is fully understood, let us consider higher rank Lie algebras. We remark that
$$
	\lm_i\leq \frac{8\pi}{\rho(A)}<4\pi,
$$
i.e. we are in a subcritical (coercive) regime. The above estimate is discussed in the next subsection \ref{subsec-spec}. Now, the uniform bound for solutions to \eqref{4.0} inside $\Omega$ has been derived in \cite{Y17} for $\textbf{A}_n, \textbf{B}_n, \textbf{C}_n$ and $\textbf{G}_2$, while $\textbf{D}_4, \textbf{F}_4$ has been discussed in \cite{KLN} for some special cases. However, as far as we are concerned just with uniform estimates in the subcritical regime and neglect the fine blow-up analysis, \cite{KLN} can be extended to the general case using the classification result for general Toda systems in \cite{LNW}. The same consideration holds for $\textbf{E}$-type Lie algebras.

\smallskip

Once we have uniform estimates inside $\Omega$, the boundary blow-up is excluded exactly as in \cite{ajy, LWZ}. This concludes the proof.

\subsection{Spectral radius of Cartan matrices} \label{subsec-spec} Finally, we discuss here $\rho(A)$ in the particular case where $A$ is a Cartan matrix associated to a symmetric simple Lie algebra.

\smallskip

The eigenvalues of $A$ are closely related to Chebyshev polynomials and have been explicitly computed in \cite{dam}. We will focus here on the $\textbf{A}_n$ case. The approach for other classical Lie algebras is similar, while the exceptional Lie algebras are finitely many and everything can be computed explicitly, see \cite{dam} for further details about eigenvalue analysis. The eigenvalues $\xi_i$ of $\textbf{A}_n$ are given by the formula
$$
	\xi_i=4\sin^2\left( \dfrac{i\pi}{2(n+1)} \right), \quad i=1,\cdots,n.
$$
It is then easy to see that
$$
	\rho(\textbf{A}_n)\in[2,4], 
$$
for any $n\geq2$. This yields the last assertion of Theorem \ref{t3} and completes its proof.

\subsection{The non-symmetric case} \label{subsec-non-symm} We show here the main steps to carry out the argument for non-symmetric Lie algebras $A$, that are
$$
\textbf{B}_{n}=\left(\begin{array}{cccccccc}
2&-1&0&0&\cdots&0\\
-1&2&-1&0&\cdots&0\\
0&-1&2&-1&\cdots&0\\
\vdots&\vdots&\vdots&\vdots&\cdots&\vdots\\
0&\cdots&\cdots&-1&2&-2\\
0&\cdots&\cdots&0&-1&2
\end{array}\right), \quad \textbf{C}_{n}=\left(\begin{array}{cccccccc}
2&-1&0&0&\cdots&0\\
-1&2&-1&0&\cdots&0\\
0&-1&2&-1&\cdots&0\\
\vdots&\vdots&\vdots&\vdots&\cdots&\vdots\\
0&\cdots&\cdots&-1&2&-1\\
0&\cdots&\cdots&0&-2&2
\end{array}\right),
$$ 

\medskip

$$
\textbf{G}_{2}=\left(\begin{array}{cccccccc}
2&-1\\
-3&2
\end{array}\right), \quad \textbf{F}_{4}=\left(\begin{array}{cccccccc}
2&-1&0&0\\
-1&2&-2&0\\
0&-1&2&-1\\
0&0&-1&2
\end{array}\right).
$$
We start by considering the following symmetric decomposition:
$$
	A=DA^s,
$$ 
where $D$ is a diagonal matrix with entries $d_i$ and $A^s$ is a symmetric matrix. We have, respectively

$$
d_i=
\begin{cases}
1, \ i=n,\\
2, \ i=1,\cdots,n-1,
\end{cases}
 \ \
\textbf{B}^s_{n}=\left(\begin{array}{cccccccc}
1&-\frac12&0&0&\cdots&0\\
-\frac12&1&-\frac12&0&\cdots&0\\
0&-\frac12&1&-\frac12&\cdots&0\\
\vdots&\vdots&\vdots&\vdots&\cdots&\vdots\\
0&\cdots&\cdots&-\frac12&1&-1\\
0&\cdots&\cdots&0&-1&2
\end{array}\right),
$$
$$
d_i=
\begin{cases}
1, \ i=1,\cdots,n-1,\\
2, \ i=n,
\end{cases}
 \  \
\textbf{C}^s_{n}=\left(\begin{array}{cccccccc}
2&-1&0&0&\cdots&0\\
-1&2&-1&0&\cdots&0\\
0&-1&2&-1&\cdots&0\\
\vdots&\vdots&\vdots&\vdots&\cdots&\vdots\\
0&\cdots&\cdots&-1&2&-1\\
0&\cdots&\cdots&0&-1&1
\end{array}\right),
$$
$$
d_i=
\begin{cases}
1, \ i=1,\\
3, \ i=2,
\end{cases}
 \  \
\textbf{G}^s_{2}=\left(\begin{array}{cccccccc}
2&-1\\
-1&\frac23
\end{array}\right),
$$
$$
d_i=
\begin{cases}
1, \ i=3,4,\\
2, \ i=1,2,
\end{cases}
 \  \
\textbf{F}^s_{4}=\left(\begin{array}{cccccccc}
1&-\frac12&0&0\\
-\frac12&1&-1&0\\
0&-1&2&-1\\
0&0&-1&2
\end{array}\right).
$$
We then also let
$$
	\lm_i^s=d_i\lm_i, \quad i=1,\dots,n.
$$

\medskip

Now, letting 
$$
A^s=(a_{ij}^s), \quad d_i v_i=u_i, \quad i=1,\dots,n,
$$ 
the system \eqref{4.0} can be rewritten as
$$
\begin{cases}
-\Delta v_{i}=\displaystyle{\sum_{j=1}^n} {a}_{ij}^s\lambda_j \frac{h_je^{d_j v_j}}{\int_{\Omega}h_je^{d_j v_j}} \ \ &\mathrm{in} \ \ \Omega,\\
v_i=0 \ \ &\mathrm{on} \ \ \partial \Omega,
\end{cases}
$$
for which the linearization leads to a problem similar to \eqref{sys1.0}, that is
$$
\begin{cases}
-\Delta\phi_i=\sum\limits_{j=1}^na_{ij}^sV_j^s\phi_j, \ \ &\mathrm{in} \ \ \Omega,\\
\phi_i=c_i\in\mathbb{R}, \ \ &\mathrm{on} \ \ \partial\Omega,\\
\int_{\Omega}V_i\phi_i=0,
\end{cases}
$$
for $i=1,\cdots,n$, where $V_j^s$ takes now the form
$$
V_j^s=\lambda_j^s\frac{h_je^{u_j}}{\int_{\Omega}h_je^{u_j}}.
$$
Since $A^s$ is a symmetric matrix we can then run the argument as in the previous subsections by using the spectral radius $\rho(A^s)$ and replacing $\lambda_j$ by $\lambda_j^s$. As for the symmetric case, to complete the proof we need to show that
$$
	\rho(A^s)\in[2,4], 
$$
for any $n\geq2$.

\

{\bf Case 1.} We consider $A^s=\textbf{B}_n^s$. Let
\begin{equation*}
\begin{split}
X_n&=|\lambda E_n-\textbf{B}_n^s|\\
&=\left|
\begin{array}{cccccccc}
\lambda-1&\frac12&0&0&\cdots&0\\
\frac12&\lambda-1&\frac12&0&\cdots&0\\
0&\frac12&\lambda-1&\frac12&\cdots&0\\
\vdots&\vdots&\vdots&\vdots&\cdots&\vdots\\
0&\cdots&\cdots&\frac12&\lambda-1&1\\
0&\cdots&\cdots&0&1&\lambda-2
\end{array}\right|.
\end{split}
\end{equation*}
We have
\begin{equation*}
X_n=(\lambda-1)X_{n-1}-\frac14X_{n-2}.
\end{equation*}
Take $\lambda\geq2$. By a straightforward computation we have
\begin{equation*}
X_n-aX_{n-1}=(\lambda-1-a)\left(X_{n-1}-aX_{n-2}\right),
\end{equation*}
with 
$$
a=\frac{\lm-1-\sqrt{(\lm-1)^2-1}}{2}>0.
$$ 
Thus, we obtain
\begin{equation*}
X_n=aX_{n-1}+(\lambda-1-a)^{n-2}\left(X_2-aX_1\right).
\end{equation*}
Now, for $\lm\geq3$ we have 
$$
X_2-aX_1>0
$$
and then, by induction we deduce $X_n>0$ for any $n\geq2$. Therefore, $\rho(\textbf{B}^s_n)<3$.

\smallskip

On the other hand, for $\lm=2$ we have 
$$
X_2-aX_1<0
$$
and then, again by induction we can get $X_n<0$ for any $n\geq2$ and hence $\rho(\textbf{B}^s_n)>2$.

\medskip

{\bf Case 2.} We consider now $A^s=\textbf{C}^s_n$. Let
\begin{equation*}
\begin{split}
X_n&=|\lambda E_n-\textbf{C}^s_n|\\
&=\left|
\begin{array}{cccccccc}
\lambda-2&1&0&0&\cdots&0\\
1&\lambda-2&1&0&\cdots&0\\
0&1&\lambda-2&1&\cdots&0\\
\vdots&\vdots&\vdots&\vdots&\cdots&\vdots\\
0&\cdots&\cdots&1&\lambda-2&1\\
0&\cdots&\cdots&0&1&\lambda-1
\end{array}\right|.
\end{split}
\end{equation*}
Take $\lambda\geq 4$. A straightforward computation yields
\begin{equation*}
X_n-aX_{n-1}=(\lambda-2-a)\left(X_{n-1}-aX_{n-2}\right),
\end{equation*}
with
$$
	a=\frac{\lambda-2-\sqrt{\lambda^2-4\lambda}}{2}>0.
$$
Thus, we obtain
\begin{equation*}
X_n=aX_{n-1}+(\lambda-2-a)^{n-2}\left(X_2-aX_1\right).
\end{equation*}
Similarly as before, one can prove that $\rho(\textbf{C}^s_n)\leq 4$ for any $n\geq 2$.

\smallskip

On the other hand, for $\lambda=2$, looking at the principal minors it is easy to see that the matrix $(\lambda E_n-\textbf{C}^s_n)$ is not positive-definite. Therefore, we conclude $\rho(\textbf{C}^s_n)>2$ for any $n\geq 2$.

\medskip

{\bf Case 3.} Finally, the spectral radius of $\textbf{G}^s_2$ and $\textbf{F}^s_4$ can be computed explicitly.

\

\begin{center}
\textbf{Acknowledgments}
\end{center}

The authors would like to thank the anonymous referee for pointing out the difference between symmetric and non-symmetric cases.

\

\

\end{document}